\documentclass[12pt,twoside]{article}
\usepackage{amsmath}
\usepackage{graphicx}
\usepackage{yhmath}
\usepackage{mathdots}
\usepackage[nottoc,numbib]{tocbibind}
\usepackage{mathtools}
\usepackage[all,cmtip]{xy}
 \usepackage{amsfonts}
\usepackage{amsthm}
\usepackage{amsmath,amscd}

\usepackage{amssymb}
\usepackage{tikz}
\usetikzlibrary{arrows.meta, decorations.markings}
\date{June 18, 2021}

\begin{document}

\centerline {\Large{\bf  On The Mutations Loops of Valued Quivers}}

\centerline{}

\centerline{}

\centerline{\bf {Ibrahim Saleh}}

\centerline{Email: salehi@uww.edu}
\centerline{University of Wisconsin Whitewater}

 \newtheorem{thm}{Theorem}[section]
 \newtheorem{cor}[thm]{Corollary}
  \newtheorem{cor and defn}[thm]{Corollary and Definition}
 \newtheorem{lem}[thm]{Lemma}
 \newtheorem{prop}[thm]{Proposition}
 \theoremstyle{definition}
 \newtheorem{defn}[thm]{Definition}
 \newtheorem{defns}[thm]{Definitions}
 \newtheorem{defns and nots}[thm]{Definitions and Notations}
 \theoremstyle{remark}
 \newtheorem{rem}[thm]{Remark}
 \newtheorem{rems}[thm]{Remarks}
 \newtheorem{exam}[thm]{Example}
 \newtheorem{exams}[thm]{Examples}
 \newtheorem{conj}[thm]{Conjecture}
 \newtheorem{que}[thm]{Question}
  \newtheorem{ques}[thm]{Questions}
  \newtheorem{ques and Conj}[thm]{Questions and Conjecture}
 \newtheorem{rem and def}[thm]{Remark and Definition}
 \newtheorem{def and rem}[thm]{Definition and Remark}
 \newtheorem{corr}[thm]{Corollary of the Proof of Proposition 5.15}
 \numberwithin{equation}{section}
\newtheorem{IbrI}{Lemma}[section]
\newtheorem{chiral}[IbrI]{Definition}
\newtheorem{IbrII}[IbrI]{Lemma}
\newcommand{\field}[1]{\mathbb{#1}}
 \newtheorem{mapProof}[thm]{Map of the Proof}
 \newcommand{\eps}{\varepsilon}
 \newcommand{\To}{\longrightarrow}
 \newcommand{\h}{\mathcal{H}}
 \newcommand{\s}{\mathcal{S}}
 \newcommand{\A}{\mathcal{A}}
 \newcommand{\J}{\mathcal{J}}
 \newcommand{\M}{\mathcal{M}}
 \newcommand{\W}{\mathcal{W}}
 \newcommand{\X}{\mathcal{X}}
 \newcommand{\BOP}{\mathbf{B}}
 \newcommand{\BH}{\mathbf{B}(\mathcal{H})}
 \newcommand{\KH}{\mathcal{K}(\mathcal{H})}
 \newcommand{\Real}{\mathbb{R}}
 \newcommand{\Complex}{\mathbb{C}}
 \newcommand{\Field}{\mathbb{F}}
 \newcommand{\RPlus}{\Real^{+}}
 \newcommand{\Polar}{\mathcal{P}_{\s}}
 \newcommand{\Poly}{\mathcal{P}(E)}
 \newcommand{\EssD}{\mathcal{D}}
 \newcommand{\Lom}{\mathcal{L}}
 \newcommand{\States}{\mathcal{T}}
 \newcommand{\abs}[1]{\left\vert#1\right\vert}
 \newcommand{\set}[1]{\left\{#1\right\}}
 \newcommand{\seq}[1]{\left<#1\right>}
 \newcommand{\norm}[1]{\left\Vert#1\right\Vert}
 \newcommand{\essnorm}[1]{\norm{#1}_{\ess}}
 \newcommand{\beq}{\begin{equation}}
\newcommand{\eeq}{\end{equation}}
\newcommand{\rarr}{\rightarrow}
\newcommand{\cA}{\mathcal{A}}
\newcommand{\cS}{\mathcal{S}}
\newcommand{\cC}{\mathcal{C}}
\newcommand{\cU}{\mathcal{U}}
\newcommand{\cR}{\mathcal{R}}
\newcommand{\RMod}{R\text{-Mod}}
\newcommand{\AMod}{A\text{-Mod}}
\newcommand{\Rep}{\text{Rep}}
\newcommand{\Aut}{\text{Aut}}
\newcommand{\XAut}{\xi\Aut}
\newcommand{\Rtzn}{R \{\theta, z, n \}}
\newcommand{\TxC}{(\Theta, \xi, \cC)}
\newcommand{\Chir}{\text{Chir}}

\tableofcontents

\begin{abstract}
 A mutation loop of a valued quiver, $Q$, is a combination of quiver automorphisms and mutations that sends $Q$  to itself. Moreover, it will be called \emph{symmetric} if it sends $Q$ to $\epsilon\sigma(Q), \epsilon \in \{-1, 1\}$ for some permutation $\sigma$. A \emph{global mutation loop} of $Q$ is a mutation loop that is symmetric for every quiver in the mutation class of $Q$. This class of relations contains all the relations of the global mutations group yield from the group action on the mutation class of $Q$. We identify which quivers have global mutation loops and provide some of them for each case.
\end{abstract}
\textbf{keywords:} Cluster Algebras,  Global Mutation Loops.

\section{Introduction}
S. Fomin and A. Zelevinsky introduced cluster algebras in [17, 8, 9, 10, 2], and a ``geometrical" version was introduced by V. Fock and A. Goncharov in [5, 6]. A cluster algebra is a  ring with distinguished sets of generators called ``\emph{clusters}". A cluster comprises two tuples of commutative variables called the $\mathcal{A}$-variables and $\mathcal{X}$-variables. Each cluster is paired with a ``valued" quiver (skew-symmetrizable matrix) to form what is called a \emph{seed}. The combinatorial operation of producing a seed from an existing one is called \emph{mutation}, which relies mainly on the quiver part of the seed. The mutation alters the quiver to produce another one. For a quiver $Q$, the set of all quivers that can be generated from $Q$ by applying sequences of mutations is called the \emph{mutation class} of $Q$, and we will denote it by $[Q]$. Mutation classes play vital roles in the cluster algebras theory and its related topics.

The transformation of $\mathcal{A}$-variables (resp. $\mathcal{X}$-variables) is called the cluster $\mathcal{A}$-transformation (resp. cluster $\mathcal{X}$-transformation).  \textit{A mutation loop} of a quiver $Q$ is a mutation sequence, formed of quiver mutations and permutations of vertices, that sends $Q$ to itself, see [12, 7]. A mutation loop defines an autonomous discrete dynamical system as the composition of cluster transformations and permutations of coordinates (equivalently permutations of vertices). The mutation loops form a cluster modular group, which acts on some geometric objects called the cluster $\mathcal{A}$- and $\mathcal{X}$-varieties and their tropicalizations by a semifield. The discrete dynamical system induced by a mutation loop occurs in these spaces. For more details, see [12], which provides one of the main motives behind the interest in learning about mutation loops.

The second and more practical reason behind the interest in learning about mutation loops comes from [13], where the action of the \emph{global mutations group}   on the classes of mutationally equivalent seeds was studied. Moreover, the authors have proven that cluster automorphisms arise naturally when those classes of ``labeled" seeds are considered as orbits under the action of the global mutations group. However, finding the group relations of the global mutation group has yet to receive the attention it deserves despite its importance in understanding the algebraic structure of the group and hence the cluster automorphisms, see [13, 14]. This paper is a step toward figuring out these relations considering the mutations classes as orbits of the group action. The main object of this paper is \emph{global mutations loops}, which are generalized relations on the global mutations group.

A \emph{symmetric mutation loop} $\mu$ of $Q$ is a  sequence of mutations such that there is a permutation $\sigma$  where $\mu(Q)=\epsilon\sigma(Q), \epsilon \in \{-1, 1\}$. One can see that such mutations sequence $\mu$  is a generalized mutation loop.  Notice that, based on the definition of mutation loops as provided in references [12] and [7], one can see that every mutation loop is actually  symmetric mutation loop with $\sigma=1$ and $\epsilon=1$. The main objective of this article is a particular class of symmetric mutation loops called \emph{global mutation loop}.

\begin{defn}
For a valued quiver $Q$, a \emph{global mutation loop} is a sequence of mutations $\mu$ that is a symmetric mutation loop for every quiver $Q'\in [Q]$. 
\end{defn}
For every single mutation $\mu_{i}$, we consider the global mutation loop $\mu^{2}_{i}$ trivial. 
The primary purpose of this article is to identify valued quivers that exhibit non-trivial global mutation loops and provide them. These valued quivers belong to the broader class of finite mutation-type quivers. The main results of this article are summarized in Theorem 1.2 below. For more details, see Example 3.2, Lemma 3.3, Proposition 3.7, and Theorem 3.10.  To facilitate a better understanding, we will introduce the necessary definitions before presenting the main results.

A quiver $Q$ is called \emph{fully cyclic} if every vertex is in a cycle. The\emph{ weight of a quiver} $Q$ is the maximum weight of its edges, and the weight of $[Q]$, denoted by $w([Q])$, is the maximum weight of the quivers in $[Q]$. Since all mutation sequences of a rank two quiver are global mutation loops, we will focus on quivers of rank three or more. 
\begin{thm} Let $Q$ be a valued quiver. Then we have the following
\begin{enumerate}
    \item If $Q$ is of rank three, then $Q$ will have global mutation loops only if  $w([Q])=4$, except the case of \begin{equation}\label{}
   \nonumber Q=
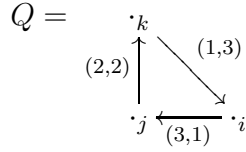
\xymatrix{\cdot_{k}\ar[dr]^{(1,3)}\\
 \cdot_{j}  \ar[u]^{(2,2)}& \cdot_{i}\ar[l]^{(3, 1)} }
  \end{equation}
  \item If $Q$ is of rank bigger than or equals to four and $w([Q])=4$, then  $Q$ has global mutation loops if and only if  each quiver $Q'\in [Q]$ is fully cyclic.
\end{enumerate}

\end{thm}
  The proofs and arguments are mainly based on  Definition 2.2 and Remarks 2.3, which is the valued quiver's version of the mutation of skew-symmetrizable matrices definition.
   In this version of the definition, mutations are applied at the vertices to alter the valued quiver by changing its edges and weights based on specific rules. Since the paper's primary focus is on quivers and their mutations, we will ignore the frozen vertices of cluster algebras.

The paper is organized as follows: we offer a concise introduction to valued quivers and their mutation. The third section provides an overview of mutation loops, including Example 3.2, Lemma 3.3, Proposition 3.7, and Theorem 3.10, which establish the category of valued quivers featuring global mutation loops. Additionally, this section includes proofs for the remaining classes that does not have non-trivial global mutation loops.

\section{Valued Quivers Mutation}

\begin{defns}
\begin{enumerate}
 
  \item \emph{An oriented valued quiver} of rank $n$ is a quadruple $Q=(Q_{0}, Q_{1}, V, d)$, where
  \begin{itemize}
 \item $Q_{0}$ is a  set of $n$ vertices labeled by $[1,n]$.
  \item  $Q_{1}$ is a set of ordered pairs of vertices, that is $Q_{1}\subset Q_{0}\times Q_{0}$ such that; $(i, i)\notin Q_{1}$ for every $i\in Q_{0}$, and if $(i, j)\in Q_{1}$, then $(j, i)\notin Q_{1}$.
  \item $V=\{(d_{ij},d_{ji})\in \mathbb{N}\times\mathbb{N} | (i,j)\in Q_{1}\}$, $V$ is called the valuation of $Q$. The weight of an edge $\alpha=(i, j)$ is the product $d_{ij}d_{ji}$ and is denoted by  $w_{i,j}$. \emph{The weight} of $Q$ is given by $w(Q)=max\{w_{ij}; (i,j) \in Q_{1}\}$.
  \item $d=(d_{1},\cdots, d_{n})$ where $d_{i}$ is a positive integer for each $i$, such that $d_{i}d_{ij}=d_{ji}d_{j}$ for every $i, j\in [1, n]$.
   \end{itemize}
  In the case of $(i,j)\in Q_{1}$, then there is an arrow oriented from $i$ to $j$, and in notation, we shall use the symbol $\xymatrix{{\cdot}_{i} \ar[r]^{(d_{ij},d_{ji})}&{\cdot}_{j}}$, or $\xymatrix{{\cdot}_{i} \ar[r]^{w_{i,j}}&{\cdot}_{j}}$ when the emphasis is on the weight of the edge. We will also use $rk(Q)$ for the rank of $Q$. The quiver $Q$ will be called \emph{simply-laced} if $w(Q)=1$. We will ignore labeling any edge of weight one. Also, a vertex $i\in Q_{0}$ is called a \emph{leaf } if there is exactly one vertex $j$ such that $w_{ij}\neq 0$ and $w_{kj}=0$, for all $k\in[1, n] \backslash \{j\}$.
  \item  A valued quiver $Q'$ is called \emph{symmetric} to $Q$ if and only if there is a permutation $\tau$ such that $Q'$ can be obtained from $Q$ by permuting the vertices of $Q$ using $\tau$ such that for every edge $\xymatrix{{\cdot}_{i} \ar[r]&{\cdot}_{j}}$ in $Q$,   the valuation $(d_{ij}, d_{ji})$ is assigned  to the edge $\xymatrix{{\cdot}_{\tau(i)} \ar[r]&{\cdot}_{\tau(j)}}$ in $Q'$.
\item A quiver $Q$ of rank $n>2$, is said to be a \emph{star quiver} if there is a vertex $i\in Q_{0}$, such that $w_{ij}\neq 0$, for $j\in [1, n]\backslash \{i\}$ and  $w_{kj}= 0$, for every $k\in [1, n]\backslash \{i\}$. In other words, each vertex is a leaf except for exactly one.
\end{enumerate}
\end{defns}
We note that every oriented valued quiver corresponds to a skew symmetrizable matrix   $B(Q)=(b_{ij})$ given by
\begin{equation}\label{}
  b_{ij}=\begin{cases} d_{ij}, & \text{ if }(i,j)\in Q_{1},\\
    0, & \text{ if }i=j,\\
-d_{ij}, & \text{ if }(j,i)\in Q_{1}.
    \end{cases}
\end{equation}
One can also see that every skew symmetrizable matrix $B$ corresponds to an oriented valued quiver $Q$ such that $B(Q)=B$.

 All our valued quivers are oriented, so in the rest of the paper, we will omit the word ``oriented". All quivers are of rank $n$ unless stated otherwise.
 We will also remove the word valued from the term ``valued quiver" when there is no confusion.

   \begin{defn}[\emph{Valued quiver mutation}]
 Let $Q$ be a valued quiver. The mutation  $\mu_{k}(Q)$ at a vertex  $k$  is defined through Fomin-Zelevinsky's mutation of the associated skew-symmetrizable matrix. The mutation of a skew symmetrizable matrix $B=(b_{ij})$ on the direction $k\in [1, n]$ is given by $\mu_{k}(B)=(b'_{ij})$, where
\begin{equation}
b'_{ij}=\begin{cases} -b_{ij}, & \text{if} \ k \in \{i,j\},\\
   b_{ij}+\text{sign}(b_{ik})\max(0, b_{ik}b_{kj}), & \text{otherwise.}
   \end{cases}
  \end{equation}

\end{defn}
The following remarks provide a set of rules that are adequate to calculate mutations of valued quivers without using their associated skew-symmetrizable matrix.
\begin{rems}
\begin{enumerate}

  \item  Let $Q=(Q_{0}, Q_{1}, V, d)$ be a valued quiver. The mutation $\mu_{k}(Q)$ at the  vertex  $k$   is described using the mutation of $B(Q)$ as follows: Let $\mu _{k}(Q)=(Q_{0}, Q'_{1}, V', d)$, we obtain $Q'_{1}$ and $V'$,   by altering  $Q_{1}$ and $V$, based on the following rules.
\begin{enumerate}
 \item replace the  pairs $(i, k)$ and $(k, j)$  with $(k, i)$ and $(j,k)$  respectively and, in the same manner, switch the components of the  ordered pairs of their valuations;
  \item if  $(i,k), (k,j)\in Q_{1}$, such that neither of $(j,i)$  or $(i,j)$ is in $Q_{1}$ (respectively $(i,j)\in Q_{1}$) add the  pair $(i, j)$ to $Q'_{1}$, and give it the valuation $(v_{ik}v_{kj},v_{ki}v_{jk})$ (respectively change its valuation to $(v_{ij}+v_{ik}v_{kj},v_{ji}+v_{ki}v_{jk})$);
 \item if $(i,k)$, $(k,j)$ and $(j, i)$ in $Q_{1}$, then we have three cases
 \begin{enumerate}
   \item if $v_{ik}v_{kj}<v_{ij}$, then keep $(j,i)$ and change its valuation to $(v_{ji}-v_{jk}v_{ki}, -v_{ij}+v_{ik}v_{kj})$;
   \item if $v_{ik}v_{kj}>v_{ij}$, then replace $(j,i)$ with $(i,j)$ and change its valuation to $(-v_{ij}+v_{ik}v_{kj}, |v_{ji}-v_{jk}v_{ki}|)$;
   \item if $v_{ik}v_{kj}=v_{ij}$,  then remove $(j,i)$ and its valuation.
 \end{enumerate}
 \item $d$ will stay the same in $\mu_{k}(Q)$.
\end{enumerate}

 \item One can see that; $\mu^{2}_{k}(Q)=Q$ and  $\mu_{k}(B(Q))=B(\mu_{k}(Q))$ at each vertex  $k\in [1, n]$ where $\mu_{k}(B(Q))$ is the mutation of the matrix $B(Q)$. For more information on mutations of skew-symmetrizable matrices, see, for example, [10, 14].
 \end{enumerate}
\end{rems}

\begin{defns}
\begin{itemize}
         \item
      The \emph{mutation class} of $Q$  is the set of all quivers that can be produced from $Q$ by applying some sequence of mutations, and it will be denoted by $[Q]$.
\item
A quiver $Q$ is called \emph{fully cyclic} if none of the vertices of $Q_{0}$  is a leaf, i.e., each vertex is a vertex in one or more cycles.  Furthermore, $[Q]$ is called \emph{fully cyclic} if every quiver in $[Q]$ is fully cyclic.

      \item  Let $\mathcal{S}_{n}$ be the symmetric group in $n$ letters. One can introduce an action of  $\mathcal{S}_{n}$ in the set of quivers of rank $n$ as follows: for a permutation $\tau$, the quiver $\tau(Q)$  is  obtained from $Q$ by permuting the vertices of $Q$ using $\tau$ such that for every edge $\xymatrix{{\cdot}_{i} \ar[r]&{\cdot}_{j}}$ in $Q$,   the valuation $(d_{ij}, d_{ji})$ is assigned  to the edge $\xymatrix{{\cdot}_{\tau(i)} \ar[r]&{\cdot}_{\tau(j)}}$ in $\tau(Q)$.
\item The \emph{weight of the mutation class} $[Q]$ is  $w([Q])=\max \{w(Q'); Q'\in [Q]\}$.

\item If $\mu=\mu_{i_{1}}\ldots \mu_{i_{k}}$ is a sequence of mutations, we will use the notations  $\{\mu\}:=\{\mu_{i_{1}}, \ldots, \mu_{k}\}$ and $[\mu_{j}]_{\mu}$ for the number of times the single mutation $\mu_{j}$ appears in $\mu$.
\end{itemize}
\end{defns}

\begin{defn}(\textbf{Cluster pattern of $[Q]$}).
The \emph{cluster pattern} $\mathbb{T}_{n}(Q)$ of the mutation class $[Q]$ is a regular $n-$ary tree whose edges are labeled by the numbers $1,2, \ldots, n$ such that the $n$ edges emanating from each vertex receive different labels. The vertices are assigned to be the elements of $[Q]$  such that the endpoints of any edge are obtained from each other by the quiver mutation in the direction of the edge label.
\end{defn}

\begin{defn}
 A quiver $Q$ is called \emph{finite mutation type} if the mutation class $[Q]$ contains finitely many quivers.
\end{defn}

\begin{defns}[Subquivers][11]
\begin{enumerate}
  \item  Let $Q$ be a quiver of rank $n$. We can obtain a sub mutation class of $[Q]$ by fixing a subset $I$ of $Q_{0}$  and apply all possible sequences of  mutations on the directions  of vertices from  the set $I$ only. In such case we say  $Q_{I}$ is  a subquiver of $Q$.  We will use $Q_{I}\leq Q$ to say $Q_{I}$ is a subquiver of $Q$. The vertices $Q_{0}\backslash I$ will be called the \emph{frozen vertices} of $Q_{I}$.

      \item
      A subquiver $Q_{I}$ is said to be of $A$-\emph{type }if the underlying graph of  $Q_{I}$  is  of $A_{n}$-type such that $w[Q_{I}]=1$.
\end{enumerate}

\end{defns}

     \begin{exams}
     \begin{enumerate}
       \item Let
     \begin{equation}\label{}
     Q=\xymatrix{ \cdot_{4}& \cdot_{3}   \ar[l]_{(2,3)}\ar[r]^{(2,3)}&   \cdot_{2}\ar[d]^{(1,2)}&\cdot_{7}\ar[l]_{(2,1)}\\
  &\cdot_{6} \ar[r]&\cdot_{1}\ar[ul]^{(6,2)}\ar[r]^{(2,3)}&\cdot_{5} }
     \end{equation}
     Consider the subquiver $Q_{I}$, with $I=\{1, 2, 3\}$,  and $d=(1, 2, 3)$. So, $rk(Q_{I})=3$ and $w(Q_{I})=12$. Here, $(6,2)$ is the valuation of the edge  $\xymatrix{{\cdot}_{1} \ar[r]&{\cdot}_{3}}$. Applying mutation at the vertex  $2$ produces the following  quiver

  \begin{equation}\label{}
      \nonumber  \mu_{2}(Q_{I})=\xymatrix{ \cdot_{4}& \cdot_{3}   \ar[l]_{(2,3)}&   \cdot_{2} \ar[l]_{(3,2)}\ar[r]^{(1,2)}&\cdot_{7}\ar[dl]^{(2, 2)}\\
  &\cdot_{6} \ar[r]&\cdot_{1}\ar[u]^{(2,1)}\ar[r]_{(2,3)}&\cdot_{5}.}
     \end{equation}

       \item Consider the quiver
       $Q=\xymatrix{ \cdot_{n}& \cdot_{n-1}\ar[l]& \cdots \ar[l]& \cdot_{2}\ar[l]&\cdot_{1}\ar[l]_{(2, 1)} }$. An easy exercise is to show that there is a sequence of mutations to move the weight two edge to any other position.

     \end{enumerate}

  \end{exams}

\begin{lem}  If $Q$ is simply-laced quiver then for every permutation $\sigma$ there is a sequence of mutations $\mu$ such that $\sigma(Q)=\mu(Q)$, for more details see [Theorem 2.6 in [1]].
\end{lem}


\section{Global Mutation loops}\label{s:2}

Fix a  quiver $Q$   with $[Q]$ as its  mutation class. Let $\mathcal{M}$ be the set of all reduced sequences (formal words) formed from the elements of the set $\{\mu_{1},\ldots, \mu_{n}\}$. Consider the action of the elements of $\mathcal{M}$  on the mutation class $[Q]$  where elements of $\mathcal{M}$ are formed by compositions of mutations. 

\begin{defns and nots}
\begin{enumerate}
\item [A.] Let $\mathcal{M}(Q)$ be the group generated by the elements of $\mathcal{M}$ subject to the relations due to its action on $[Q]$. The group $\mathcal{M}(Q)$ is called \emph{the global mutation group} of $Q$. 

 \item [B.] Consider the following particular types of elements of $\mathcal{M}(Q)$.
\begin{enumerate}

 \item [1.] A sequence of mutations  $\mu$ is called \emph{homogeneous } if whenever  $\mu=\mu^{(1)}\mu^{(2)}$, such that $\mu^{(i)}\neq 1$ for $i=1,2$, then we have
   $\{\mu^{(1)}\} \cap \{\mu^{(2)}\} \neq \emptyset$. Otherwise $\mu$ will be called \emph{inhomogeneous} and $\mu^{(i)}$ will be called \emph{isolated subsequence} of mutations of $\mu$ for $i=1,2$.
    \item [2.] A sequence of mutations  $\mu$ will be \emph{full}  if $\{\mu\}=\{\mu_{1},\ldots, \mu_{n}\}$.

 \item [3.] Let $\mu$ be a mutation sequence such that $\{\mu \}$ contains more than one element. Then $\mu$ is called \emph{symmetric mutation loop} of $Q$ if $\mu(Q)=\epsilon\tau(Q), \epsilon\in \{-1, 1\}$ for some permutation $\tau$. The subset of $\mathcal{M}(Q)$ of all symmetric mutation loops is called the \emph{symmetric mutations set} of $Q$, and it will be denoted by $\mathbf{L}_{Q,\mathcal{S}}$.
 \item [4.] The set of \emph{global symmetric loops} on $Q$ is denoted by  $\mathbf{B}_{[Q],\mathcal{S}}$ and it will be defined as follows
 \begin{equation}\label{}
  \nonumber \mathbf{B}_{[Q],\mathcal{S}}=\bigcap_{Q'\in [Q]} \mathbf{L}_{Q',\mathcal{S}}.
 \end{equation}
\end{enumerate}
\end{enumerate}
  \end{defns and nots}
We will omit the word symmetric from the term  ``global symmetric loops" as long as there is no confusion. The rest of the paper aims to identify quivers with non-trivial global loops and provide an explanatory statement for the remaining cases.

\begin{exams} [Quivers with global loops] In the following, we provide a list of quivers with global loops based on their ranks.
\begin{enumerate}
\item Let $Q$ be the  quiver $\xymatrix{\cdot_{i} \ar[r]^{(m, n)}& \cdot_{j}}$. Then we have $\mathbf{B}_{[Q], \mathcal{S}}=\{1, \mu_{i}\}=\{1, \mu_{j}\}$.
  \item If $Q$ is a  periodic quiver, i.e., $[Q]=\{\pm Q\}$, of finite mutation type of rank three, then every mutation sequence is a global loop.
  \item  Let $Q$ be the  quiver

 \begin{equation}
 \nonumber   
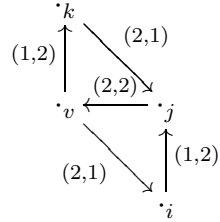
\xymatrix{\cdot_{k}\ar@{-}[dr]\ar[dr]^{(2, 1)}\\
 \cdot_{v}  \ar[u]^{(1, 2)}\ar[dr]_{(2, 1)}& \cdot_{j}\ar[l]_{(2,2)} \\
 & \ar[u]_{(1, 2)}\cdot_{i}} \
\end{equation}
Since the class $[Q]$ is very small, containing only four quivers, it would be a good exercise  to show that   $\mu_{k}\mu_{i}, \mu_{i}\mu_{k}, \mu_{v}\mu_{j}$ and $\mu_{j}\mu_{v}$ are the non-trivial global loops of $[Q]$.

\item Let $Q$ be the quiver
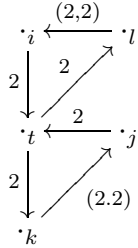
\begin{equation}\label{}
  \nonumber \xymatrix{\cdot_{i}\ar[d]_{2}&  \ar[l]_{(2,2)} \cdot_{l} \\
\cdot_{t}\ar[ur]^{2} \ar[d]_{2} & \ar[l]_{2}\cdot_{j} \\
	\cdot_{k}  \ar[ur]_{(2. 2)}}  \
\end{equation}

Since the mutation class $[Q]=\{Q\}$, one can see that $\mathbf{B}_{[Q], \sigma}$ is generated by $\{\mu_{i}, \mu_{l}, \mu_{t}, \mu_{j}, \mu_{k}\}$.
\end{enumerate}
\end{exams}

In the following lemma, we identify all rank three quivers with global loops.
\begin{lem} Let $Q$ be a quiver  of rank 3. Then we have the following.
\begin{enumerate}
  \item If $w([Q])=m$, where $ m=1,$ or $ 2$, then  $\mathcal{M}(Q)$ does not  have any global loops.
  \item If $w([Q])=4$, then $\mathcal{M}(Q)$ contains global loops except the case  \begin{equation}\label{}
   \nonumber Q=
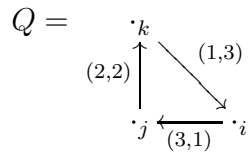
\xymatrix{\cdot_{k}\ar[dr]^{(1,3)}\\
 \cdot_{j}  \ar[u]^{(2,2)}& \cdot_{i}\ar[l]^{(3, 1)} }
  \end{equation}.
\end{enumerate}
\end{lem}

\begin{proof}
\begin{enumerate}
  \item Assume that  $w[Q]=m; \ \text{where} \ m=1$, or $ 2$. Without loss of generality,  let $Q$ be  a cyclic   quiver  with the following underlying graph
  \begin{equation}\label{}
    \xymatrix{\cdot_{k}\ar[dr]^{m}\\
 \cdot_{i}  \ar[u]& \cdot_{j}\ar[l]^{m} }
  \end{equation}
   \begin{table}
\caption{Counter examples of quivers for rank 3 with weight $m=1$ or $2$}
\begin{tabular}{ l l }
    \hline
    \bf $\mu $ & \bf Counter examples of quivers\\
    \hline
    $\mu_{i}\mu_{j}\mu_{k}\mu_{i} $ &  $\xymatrix{\cdot_{i}&\ar[l]  \cdot_{k} \ar[r]^{m} & \cdot_{j} } $ \\
    $\mu_{i}\mu_{j}\mu_{k}\mu_{i}$ & $\xymatrix{\cdot_{i} \ar[r]&  \cdot_{k} \ar[r]^{m} & \cdot_{j} } $\\
    $\mu_{j} $ & $\xymatrix{\cdot_{i}&\ar[l]_{m}  \cdot_{k} \ar[r] & \cdot_{j} } $\\
        \hline
\end{tabular}
\end{table}
 One can see that if $\mu \in \mathbf{B}_{Q, \mathcal{S}}$ that is not an involution relation,  then $\mu$ must break and rebuild the cyclic quiver or only reverse of the arrows. So, $\mu$ must be one of the mutation sequences in the left-side column in Table 1.1 or formed of one or more of them. Moreover, a short calculation shows that the right-side quivers in Table 1.1 provide counter-examples of quivers where $\mu$ does not work as a global relation in $[Q]$. The case of $3$ we have $[Q]$ is actually empty assuming it is finite.
  \item Table 1.2 below provides proof of the existence of global loops for all five cases of $Q$ that have global mutation loops. Finally, let
  \begin{equation}\label{}
   \nonumber Q= \xymatrix{\cdot_{k}\ar[dr]^{(1,3)}\\
 \cdot_{j}  \ar[u]^{(2,2)}& \cdot_{i}\ar[l]^{(3, 1)} }.
  \end{equation}

  One can see that $Q$ is mutationally equivalent to the quiver in (3.1) with $m=3$. Then from Part 1 above, we have $Q$ has no global mutation loops.
  \begin{table}
\caption{Global loops for rank 3 quivers of weight 4 }
  \begin{tabular}{ l l }
    \hline
    \bf Mutation classes of weight 4 & \  \  \  \bf Global mutation loops\\
    \hline
    $\xymatrix{\cdot_{k}\ar[dr]^{(2,2)}\\
 \cdot_{i}  \ar[u]^{(2,2)} & \cdot_{j}\ar[l]^{(2,2)} } $ & \  \  \  \ \ $\mathcal{M}(Q) $\\

    $\xymatrix{\cdot_{k}\ar[dr]^{(1,4)}\\
 \cdot_{i}  \ar[u]^{(2,1)}& \cdot_{j}\ar[l]^{(2,1)} }$ & \  \  \ \  \ $\mu_{ijki}, \mu_{ikji} $\\

   $\xymatrix{\cdot_{k}\ar[dr]^{(2,1)}\\
 \cdot_{j}  \ar[u]^{(2,2)}& \cdot_{i}\ar[l]^{(1,2)} } $ & \  \  \  \  $\mu_{ijki}, \mu_{ikji} $\\

   $\xymatrix{\cdot_{k}\ar[dr]^{(2,2)}\\
 \cdot_{i}  \ar[u]^{(4,1)}& \cdot_{j}\ar[l]^{(1,4)} } $ & \  \  \  \ $\mathcal{M}(Q) $\\
        \hline
\end{tabular}
\end{table}

\end{enumerate}
\end{proof}

The quivers referred to in the following Lemma are denoted as \emph{leading quivers}. The purpose of this lemma is to establish relationships between $w([Q])$ and the structure of certain quivers within $[Q]$.

\begin{lem} $[15, \text{Lemma 3.11}]$. Let $Q$ be a  mutationally finite quiver. Then the weight of $[Q]$ is determined as follows
\begin{enumerate}
  \item  $w[Q]=2$  if and only if  $rk(Q)=2$ or $[Q]$ contains a quiver  that has exactly one  edge, say $(i, j)$, of weight $2$,   such that the sub-quiver    $Q_{\{Q_{0} \setminus \{i, j\}\}}$ is of $A$-type or in otherworld, $Q$ is mutation-equivalent to a quiver of the type $B_{n}, C_{n}$ or $F_{4}$.
  \item $w[Q]=3$ if and only if $Q$  is the  quiver $ \xymatrix{\cdot_{i} \ar[r]^{(3, 1)} & \cdot_{j}}.$
  \item    If $rk(Q)>3$ then every quiver $Q'\in [Q]$  of weight 4  satisfies that  edges of weight  4 in $Q'$  appear in  a cyclic subquiver that is symmetric to  one of the following

 \begin{equation}
\nonumber  (a)  \  \   \ Q_{a, x}: \  \xymatrix{
\cdot_{v} \ar[d]_{(x, 1)} & \ar[l]_{(1, x)}\cdot_{j} \\
	\cdot_{k}  \ar[ur]_{(2. 2)}} \  \  \  \ \ \ \  \
 (b)  \  \  \ Q_{a}: \  \xymatrix{
\cdot_{v} \ar[d]_{(1, 2)} & \ar[l]_{(1, 2)}\cdot_{j} \\
	\cdot_{k}  \ar[ur]_{(4. 1)}} \  \  \  \
 \end{equation}
 \begin{equation}\label{}
 (c) \  Q_{c, t}: \  \xymatrix{\cdot_{k}\ar[dr]^{(t, 1)}\\
 \cdot_{v}  \ar[u]^{(1, t)}\ar[dr]_{(1, 2)}& \cdot_{j}\ar[l]_{(2, 2)} \\
 & \ar[u]_{(2, 1)}\cdot_{l} }\\ \  \  \ \  \  \
 (d)  \  Q_{d}:  \  \  \  \ \  \ \xymatrix{\cdot_{v}\ar[dr]\\
 \cdot_{i}  \ar[u]\ar[dr]_{(1, 3)}& \cdot_{j}\ar[l]_{(2, 2)}. \\
& \ar[u]_{(3, 1)}\cdot_{l} } \\
 \end{equation}

Where $x=1, 2, 3$ or $4$  and  $t=1$ or $ 2$, such that edges of weight four are not connected outside their cycles.

 \item  (Quivers with more than one edge of weight 4). If $Q'$ has more than one edge of weight 4, then it is symmetric to one of the following cases.
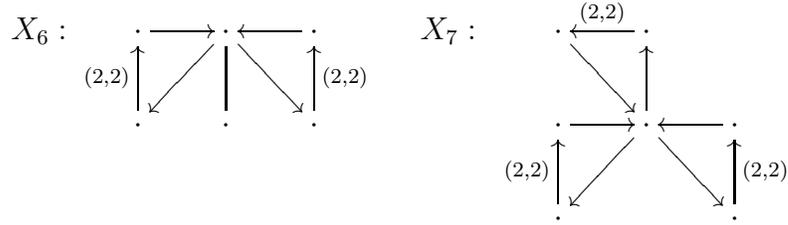
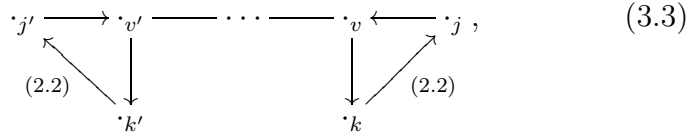
\begin{enumerate}
  \item   \begin{equation}\label{}
               \nonumber X_{6}: \xymatrix{
\cdot \ar[r] &\cdot\ar[dl] \ar[dr] & \ar[l]\cdot \\
\cdot\ar[u]^{(2,2)}&\cdot\ar@{-}[u]&\cdot  \ar[u]_{(2,2)}}
            \ \ \ \  X_{7}:
              \  \xymatrix{\cdot\ar[dr]&\cdot\ar[l]_{(2,2)}\\
\cdot \ar[r] &\cdot\ar[dl] \ar[dr]\ar[u] & \ar[l]\cdot \\
\cdot\ar[u]^{(2,2)}&&\cdot  \ar[u]_{(2,2)}}
            \end{equation}
  \item $Q^{l}_{a, 1}$:

            \begin{equation}\label{}
             \  \xymatrix{
\cdot_{j'}\ar[r]&\cdot_{v'} \ar[d]\ar@{-}[r] &\cdots\ar@{-}[r]&\cdot_{v} \ar[d] & \ar[l]\cdot_{j} \\
&\cdot_{k'}\ar[ul]^{(2. 2)}&&\cdot_{k}  \ar[ur]_{(2. 2)}},
            \end{equation}
where the subquiver connecting $v$ and $v'$ is of $A$-type of rank $l\geq1$.

\item
   $Q^{l}_{a}$ or $Q^{l}_{a, 2}$ where the weight $4$ is of valuation  $(4, 1)$ or $(2, 2)$ respectively with the obvious changes

            \begin{equation}\label{}
             \  
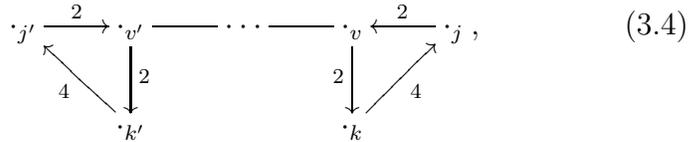
\xymatrix{
\cdot_{j'}\ar[r]^{2}&\cdot_{v'} \ar[d]^{2}\ar@{-}[r] &\cdots\ar@{-}[r]&\cdot_{v} \ar[d]_{2} & \ar[l]_{2}\cdot_{j} \\
&\cdot_{k'}\ar[ul]^{4}&&\cdot_{k}  \ar[ur]_{4}},
            \end{equation}
           also the subquiver connecting $v$ and $v'$ is of $A$-type of rank $l\geq 1$.

   \end{enumerate}

\end{enumerate}

\end{lem}

\begin{lem} Let  $Q$ be a quiver with $w([Q])\leq 2$ and $rk(Q)\geq3$. Then every global loop of $Q$ is full. In particular, the statement is true if $\mathcal{S}_{n-1}\subset \mathcal{M}(Q)$.
\end{lem}
\begin{proof} The proof is divided into two parts based on $w([Q])$.
\begin{enumerate}
  \item Assume that $w([Q])=1$. Then $\mathcal{S}_{n}\subset \mathcal{M}(Q)$. Let  $\mu\in \mathbf{B}_{Q', \mathcal{S}}$ for every $Q'\in [Q]$. Since $\mu$ is not an involution, then $\{\mu\}$ contains at least two elements. The case of $\{\mu \}$ containing exactly two elements is obviously can not be a global loop as $rk(Q)\geq3$, as counter quivers can be found easily. Now assume that $\{\mu\}$ contains three or more elements. Now, let  $i, j$ and $k$ be vertices in  $\{\mu\}$. Since $\mathcal{S}_{n}\subset \mathcal{M}(Q)$, then we can  identify a quiver $Q'\in [Q]$ such that $Q''=\sigma (Q')$, for some  permutation 
$\sigma$. Where $\sigma$ rearranges $Q'$ so that  the vertices $i, j$ and $k$ form a subquiver $Q^{\star}$ of $Q''$.  The subquiver $Q^{\star}$ is isolated from the rest of  $Q'$ by a vertex $t$ where  $t\notin \{\mu\}$. Here is how $Q^\star$ would appear as a subquiver of $Q''$
  \begin{equation}\label{}
    \xymatrix{\cdot_{i}\ar[r]& \cdot_{j} \ar[r] &\cdot_{k}\ar[r] &\cdot_{t}\ar[r]&\cdots }.
  \end{equation}

One can see that the quiver $Q''$ is not going to be altered by any mutation sequence that does not contain $\mu_{i}, \mu_{j}$ or $\mu_{k}$. Therefore, if $\mu$ is a global loop, it must contain a sub sequence of mutations that is a  global loop on the subquiver $Q''$ which is impossible, thanks to Part 1 of Lemma 3.3.
\item Assume that $w([Q])=2$. Then $[Q]$ contains a quiver $Q'$ with an underlying graph of $A_{n}$-type with exactly one single edge of weight 2, thanks to  Part 1 of Lemma 3.4. One can see that $Q'$ has the following two features

\begin{enumerate}
  \item The edge of weight two  can move to any position on $Q'$, thanks to Part 2 of Example 2.8;
  \item Without loss of generality, if we fix the edge of weight two as a leaf, say  $\xymatrix{\cdot_{n-1}\ar[r]^{2}&\cdot_{n} }$. Then the vertices in $[1, n-1]$ are permutable, i.e., $\mathcal{S}_{n-1}\subset \mathcal{M}(Q)$.
\end{enumerate}
Now, if $rk(Q')= 3$, then $Q$ does not have any global loops, thanks to Part 1 of Lemma 3.3. Assume that  $rk(Q) \geq 4$. Suppose that $\mu$ is not a full global loop. Then  at least one vertex, $t$, that is not represented in $\{\mu\}$. If $t=n$, move the edge of weight 2 to the other end so $t$ becomes free to move anywhere in the $A_{n-1}$-type subquiver of $Q'$. So whether $t$ is a leaf or not, we can move it at the position of any other vertex. Therefore, one can isolate a subquiver of three vertices, such as in (3.5),  from the rest of the quiver by the vertex  $t$ using the fact that $\mathcal{S}_{n-1}\subset \mathcal{M}(Q)$. Then, using the same argument as Part 1 of this lemma, we conclude that such a vertex $t$ does not exist. Then $\mu$ must be full.

\end{enumerate}
\end{proof}
The subsequent lemma is adapted from [15] to better suit the context of this article.
\begin{lem}[15, Lemma 3.6] A quiver $Q$ is of infinite mutation type if and only if there exists a quiver $Q'\in [Q]$ which contains a $3$-cycle $Q^{\star} \leq Q$, where
 \begin{equation}\label{}
   \nonumber Q^{\star}= \xymatrix{\cdot_{i}\ar[dr]^{w_{1}}&\\
\cdot_{k}\ar[u]^{w_{2}}& \cdot_{j} \ar[l]^{w_{3}} & },
 \end{equation}
 such that $w(\mu_{x}\mu(Q^{\star}))>w(\mu(Q^{\star}))$ for every sequence of mutations $\mu$ and for  $x \in \{i, j, k\}$.
\end{lem}

\begin{prop}
Let $Q$ be a quiver. Then we have the following.
\begin{enumerate}
  \item If $Q$ is of infinite mutation type, then $\mathcal{M}(Q)$ does not have any non-trivial global loops.
   \item If $Q$ is a   star quiver, then $\mathcal{M}(Q)$ does not have any non-trivial global loop.
  \end{enumerate}
\end{prop}

\begin{proof}
\begin{enumerate}
  \item Suppose that $Q$ is of infinite mutation type. Then $[Q]$ contains a quiver $Q'$ with an unbounded, cyclic subquiver. Hence, from Lemma 3.6, there is a  mutation sequence $\mu$ such that $w(\mu_{i}\mu(Q'))>w(\mu(Q))$ for every $i \in [1, n]$. Then $[Q]$ contains a quiver with weights that are impossible to be reproduced using any sequence of mutations. Therefor, $\mathbf{B}_{[Q], \mathcal{S}}$ can not  contain any non-involutions loops.
  \item The proof of this case is provided for a rank four star quiver, and the higher rank is the same argument. Let $Q$ be a star quiver of rank 4. Then $[Q]$ contains a quiver $Q'$, which has a subquiver with the following underlying graph

  \begin{equation}\label{}
  \nonumber 
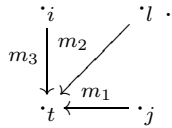
\xymatrix{\cdot_{i}\ar[d]_{m_{3}}&  \ar[dl]_{m_{2}} \cdot_{l} \\
\cdot_{t} & \ar[l]_{m_{1}}\cdot_{j}}.
\end{equation}
To finish the proof, Table 1.3  below provides the list of possible mutation loops on $Q$ and a counter example of quivers in $[Q]$ for each case.
    \begin{table}
\caption{Counter examples of quivers for rank 4 star quiver }
\begin{tabular}{ l l }
    \hline
    \bf $\mu $ & \bf Counter examples\\
    \hline
    $\mu_{t} $ &  $\xymatrix{\cdot_{i}\ar[d]_{m_{3}}&   \cdot_{l} \\
\cdot_{t}\ar[ur]^{m_{2}} & \ar[l]_{m_{1}}\cdot_{j}}$ \\
    $\mu_{t}^m(\mu_{i}\mu_{j}\mu_{l})^k, (\mu_{i}\mu_{j}\mu_{l})^k\mu_{t}^m,     
     (\mu_{i}\mu_{j}\mu_{l})^k$ \ & $\xymatrix{\cdot_{i}\ar[r]^{m_{5}}&\cdot_{j}\ar[dl] _{m_{1}}\\
\cdot_{t} \ar[u]^{m_{3}}\ar[r]_{m_{2}} & \cdot_{l}\ar[u]_{m_{4}}}$\\
            \hline
\end{tabular}

Where $m, m_x; x=1, \cdots, 5$, and $k$ are non-negative integers.
\end{table}

\end{enumerate}
\end{proof}

 \begin{lem}
If $\mathcal{S}_{n}\subset \mathcal{M}(Q)$, then every global loop is homogenous, i.e., every global loop does not contain any isolated subsequences.
\end{lem}
\begin{proof}

              We have two cases: The first is when $Q$ has a subquiver of $A_{4}$-type connected to $Q$ at exactly one vertex. Let $\mu=\mu_{i_{1}}\cdots \mu_{i_{k}}$ be a global loop on $[Q]$  with an isolated subsequence $\mu^*$  including $\mu_{1}, \mu_{2}, \mu_{3}$ and  $\mu_{4}$. Without loss of generality, assume that $\mu^*$ appears at the end of $\mu$. Recall that $\mathcal{S}_{n}\subset \mathcal{M}(Q)$, which makes it possible to use mutations to move and rearrange the vertices and customize the shape of the underlying graph to be in the following shape.

          \begin{equation}\label{}
  \nonumber 
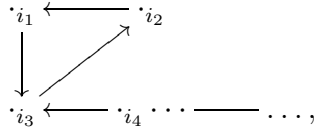
\xymatrix{\cdot_{i_{1}}\ar[d]&  \ar[l] \cdot_{i_{2}} \\
\cdot_{i_{3}}\ar[ur]  & \ar[l]\cdot_{i_{4}}\cdots\ar@{-}[r]& \ldots,}
\end{equation}
 One can see that $\mu$ will break and rebuild the triangular cycle. Then the possible choices for the action of $\mu^*$ on the subquiver of the vertices $\{i_{1}, i_{2}, i_{3}, i_{4}\}$ are equivalent to  the forms: $\mu_{i_{1}}\mu_{i_{3}}\mu_{i_{2}}\mu_{i_{1}}$,  $\mu_{i_{1}}\mu_{i_{2}}\mu_{i_{3}}\mu_{i_{1}}$, $\mu_{i_{2}}\mu_{i_{3}}\mu_{i_{1}}\mu_{i_{2}}, \mu_{i_{2}}\mu_{i_{1}}\mu_{i_{3}}\mu_{i_{2}}$ $\mu_{i_{3}}\mu_{i_{2}}\mu_{i_{1}}\mu_{i_{3}}$  and $\mu_{i_{3}}\mu_{i_{1}}\mu_{i_{2}}\mu_{i_{3}}$. Now we will provide a counter-example of a quiver from $[Q]$ for each case. Consider the following quiver from $[Q]$,

     \begin{equation}\label{}
  \nonumber \xymatrix{\cdot_{i_{z}}& \ar[l] \ar[r] \cdot_{i_{x}}&\cdot_{i_{y}}& \ar[l]\cdot_{i_{4}}\cdots\ar@{-}[r]& \ldots}.
 \end{equation}
If $z=1$ then $x=2$ and $y=3$ for $\mu_{i_{1}}\mu_{i_{3}}\mu_{i_{2}}\mu_{i_{1}}$ and $x=3$ and $y=2$ for $\mu_{i_{1}}\mu_{i_{2}}\mu_{i_{3}}\mu_{i_{1}}$. And if $z=2$ then $x=1$ and $y=3$ for $\mu_{i_{1}}\mu_{i_{3}}\mu_{i_{2}}\mu_{i_{1}}$ and $x=3$ and $y=1$ for $\mu_{i_{1}}\mu_{i_{2}}\mu_{i_{3}}\mu_{i_{1}}$. The case of $\mu_{i_{3}}\mu_{i_{2}}\mu_{i_{1}}\mu_{i_{3}}$  and $\mu_{i_{3}}\mu_{i_{1}}\mu_{i_{2}}\mu_{i_{3}}$ is similar.

The second case is when $Q$ is of $E_{6}$-type. One can see that the following quiver, $Q'$,  is in $[Q]$
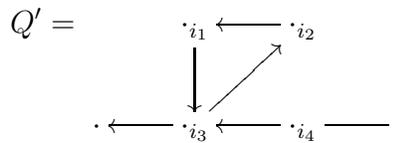
  \begin{equation}\label{}
  \nonumber Q'=\xymatrix{&\cdot_{i_{1}}\ar[d]&  \ar[l] \cdot_{i_{2}} \\
\cdot& \ar[l]\cdot_{i_{3}}\ar[ur]  & \ar[l]\cdot_{i_{4}}\ar@{-}[r]& \cdot}
\end{equation}
We use the same argument we used in the first case above,  with the quiver $Q'$,  where the counter examples of quivers are encoded in the following quiver.

 \begin{equation}\label{}
  \nonumber Q'=
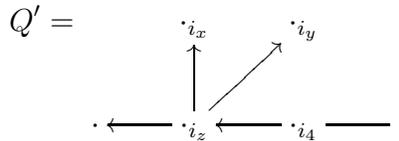
\xymatrix{&\cdot_{i_{x}}&  \cdot_{i_{y}} \\
\cdot& \ar[l]\cdot_{i_{z}}\ar[ur] \ar[u] & \ar[l]\cdot_{i_{4}}\ar@{-}[r]& \cdot}
\end{equation}
with $x, y, \ \text{and} \ z \in \{1, 2, 3\}$, which finishes the proof.
\end{proof}

\begin{rem} If $Q$ is of finite mutation type quiver, then   $[Q]$   is fully cyclic if and only if $[Q]$ contains one of the following quivers
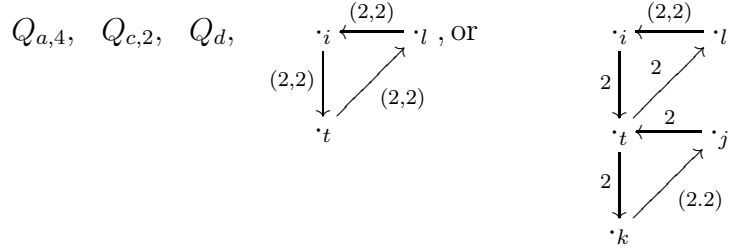
\begin{equation}\label{}
  \nonumber  Q_{a, 4}, \ \  Q_{c, 2}, \ \ Q_{d}, \ \
  \xymatrix{\cdot_{i}\ar[d]_{(2, 2)}&  \ar[l] _{(2, 2)}\cdot_{l} \\
\cdot_{t}\ar[ur]_{(2, 2)} },\text{or}   \ \ \ \ \ \ \ \ \ \ \    \xymatrix{\cdot_{i}\ar[d]_{2}&  \ar[l]_{(2,2)} \cdot_{l} \\
\cdot_{t}\ar[ur]^{2} \ar[d]_{2} & \ar[l]_{2}\cdot_{j} \\
	\cdot_{k}  \ar[ur]_{(2. 2)}}
\end{equation}
\end{rem}

In the following, we will provide the main statement of this article.

\begin{thm}  Let $Q$ be a valued quiver. Then we have the following.
\begin{enumerate}
\item If $w([Q])\leq 2$ and $rk(Q)\geq 3$, then $[Q]$  has no global loops.
\item If $w([Q])=3$,  then every sequence of mutations is a global loop.
  \item If $w([Q])=4$ and $rk(Q)>3$, then $\mathcal{M}(Q)$  has  global loops if and only if $[Q]$ is fully cyclic. In such case, $[Q]$ contains only two forms of underlying graphs.
\end{enumerate}
\end{thm}
\begin{proof}
\begin{enumerate}
  \item
  \begin{enumerate}
    \item 
    
  The case of  $rk(Q)=3$ is proved in  Part 1 of Lemma 3.3. The rest of the proof is divided into two cases based on $w([Q])$.
         Assume that $w([Q])=1$.  Then $\mathcal{S}_{n}\subseteq \mathcal{M}(Q)$. Let $\mu$ be a global loop of $[Q]$. Then $\mu$ is full and homogenous, thanks to Lemma 3.5 and Lemma 3.8. Assume that $\mu=\mu^{\star}_{2}\mu_{1}^{\star}$,
          where $\{\mu\}=Q_{0}$ and each of $\mu^{\star}_{1}$ and $\mu^{\star}_{2}$ is a subsequence of mutations of $\mu$.  In particular assume that  $\{\mu_{1}^{\star}\}=\{\mu_{i_{1}}, \mu_{i_{2}}, \mu_{i_{3}}, \mu_{i_{4}}\}$. We have the following remarks on $\mu$
  \begin{enumerate}
   \item  the subsequence $\mu^{\star}_{2}$  could be identity, depending on the rank of $Q$;
    \item  in case that  $\{\mu^{\star}_{2}\}$ is not empty, since the subsequence $\mu_{1}^{\star}$ is not isolated then $\{\mu^{\star}_{1}\}\cap \{\mu^{\star}_{2}\}$ is not empty;
       \item if $\{\mu^{\star}_{2}\}$ is not empty, then there are at least two vertices, say $z \in \{\mu^{\star}_{2}\}$ and $t \in \{\mu^{\star}_{1}\}$ such that $[\mu_{t}]_{\mu}\neq [\mu_{z}]_{\mu}$. To show that, one can use mutations and  $S_{n}\subset \mathcal{M}(Q)$ to customize a quiver $Q' \in [Q]$ with the following features
    \begin{itemize}
      \item $Q'$ has an underlying graph of $A_{n}$-type with two 3-cycles subquivers; one 3-cycle at each end;
      \item One of the 3-cycle subquivers contains both  $z$ and $t$ as vertices;
      \item $\mu_{t}$  appears in $\mu^{\star}_{2}$ before $\mu_{z}$.
    \end{itemize}
     So,  $\mu_{t}$ is in the position where applying it breaks its 3-cycle, and it has already been applied within $\mu^{\star}_{1}$. Furthermore, since $\mu$ is redundant-free then $[\mu_{t}]_{\mu}> [\mu_{z}]_{\mu}$.
  \end{enumerate}

 Since $Q$ is a simply-laced finite mutation type quiver, then one can find a quiver
$Q''$ in $[Q]$ that has the shape
  \begin{equation}\label{}
\xymatrix{\cdot_{i_{1}}\ar[d]&  \ar[l] \cdot_{i_{2}} \\
\cdot_{i_{3}}\ar[ur]  & \ar[l]\cdot_{i_{4}}\cdots\ar@{-}[r]& \ldots.}
\end{equation}
There are three possible ways for $\mu^{\star}_{1}$ to can act on the 3-cycle $[t_{1}, t_{2}, t_{3}]$
\begin{enumerate}
\item [Case 1.]  break and rebuild the 3-cycle  completely such that none of the vertices travel outside the 3-cycle over the path of $\mu$, i.e., $\mu(Q'')$ has a subquiver that is symmetric to the subquiver (3.6) with vertices  $i_{1}, i_{2},$ and  $i_{3}$.
\item [Case 2.] break and rebuild the 3-cycle  completely such that some of the vertices travel outside the 3-cycle and get back over the path of $\mu$, i.e., at least one vertex of the 3-cycle, $[t_{1}, t_{2}, t_{3}]$, is replaced over the path of  $\mu(Q'')$.
\item [Case 3.]  break the 3-cycle and then rebuild it using other vertices, i.e., the 3-cycle appears in   $\mu(Q'')$ but in the other end where all vertices $i_{1}, i_{2},$ and  $i_{3}$ are replaced.
\end{enumerate}

\begin{itemize}

\item Case 1: One can see that the 3-cycle (3.6) vertices will not leave it over the path of $\mu$. Then $\mu$ will behave as a global loop on the subquiver (3.6). However, simply laced 3-cycles have no global loops, thanks to Part 1 of  Lemma 3.3.

\item Case 2:  Suppose that $i_{2}$ travels outside the 3-cycle and swaps position  with another vertex, say $t$, and then it returns to the 3-cycle at the end over the path of $\mu$ on $Q$. Fix $i_{j}\neq i_{2}$ such that $\mu_{i_{j}}$ to be the last single mutation to appear in $\mu$.
 Use mutations and the fact that $S_{n}\subset \mathcal{M}(Q)$  to customize a quiver  $Q'$  in $[Q]$ such that it has the same underlying graph of $Q$ with $t$ swapped vertices with $i_{j}, j\neq 2$, where $\mu_{i_{j}}$ is the last single mutation that appears in $\mu$. Here, the 3-cycle $[t_{1}, t_{2}, t_{3}]$ is formed right after applying $\mu_{i_{j}}$ on the path of $\mu(Q)$. One can see that rebuilding the 3-cycle $[t_{1}, t_{2}, t_{3}]$ in $Q'$ using $\mu$ will be interrupted and would not be completed  which means $\mu(Q')$ is not symmetric to $Q'$. Then such a global loop $\mu$ does not exist.

  \item Case 3: Without loss of generality assume that $\mu^{\star}_{2}$ is not trivile. In this case, the action of $\mu$ causes all the 3-cycles, $[t_{1}, t_{2}, t_{3}]$, vertices to be swapped with other vertices or the 3-cycle that appear in the other end. Without loss of generality, suppose that the new 3-cycle is formed of vertices from $\{\mu^{\star}_{2}\}$.
       Since $\mu$ is a homogenous mutation loop, then $\{\mu^{\star}_{2}\} \cap \{\mu^{\star}_{1}\}$ is non-empty. Now, let $\mu_{t}$ be the last single mutation to appear in $\mu$ from  $ \{\mu^{\star}_{2}\}\cap \{\mu^{\star}_{1}\}$.
        Use mutations and permutations to customize a quiver  $Q''$  in $[Q]$ such that it has the same underlying graph of $Q$ with the vertex $t$  is swapped with $z$  and $t_{1}$ becomes a leaf in $Q''$. On the path of $\mu$ on $Q$, the 3-cycle will be rebuilt right after applying $\mu_{t}$  and it will stay with no change since it is the last single mutation to appear in $\{\mu^{\star}_{2}\} \cap \{\mu^{\star}_{1}\}$. However, on the other hand, the 3-cycle would not be built in $\mu(Q'')$ because over the path of $\mu$ on $Q''$, the vertex $t$ is exchanged with $z$ and  $[\mu_{z}]_{\mu}\neq [\mu_{t}]_{\mu}$. Hence, $\mu(Q')$ is not symmetric to $Q''$.

\end{itemize}

          \item Assume that $w([Q])=2$. Then if there is any global loop, it must be full, thanks to  Lemma 3.5. The mutation class of $[Q]$ contains a quiver $Q'$ with an $A_{n}$-type underlying graph that has one single edge of weight 2, say $\xymatrix{\cdot_{i} \ar[r]^{2}&  \cdot_{j} }$.  Without loss of generality, assume that the edge of  weight two in $Q'$ is a leaf. Using the fact that the weight two edge can move to any spot and $S_{n-1}\subset \mathcal{M}(Q)$, we can imitate a similar proof of the simply-laced case.

  \end{enumerate}
 \item Assume that $w([Q])=3$. Then Part 2 of Lemma 3.4 guarantees that  $Q$ has only two vertices. Therefore, every sequence of mutations of $Q$ is a global loop.
\item The proof of this case is a case-by-case process based on the \emph{leading quivers}, (the quivers in Lemma 3.4). ``$\Leftarrow$". This case is proved by Table 1.4.

  \begin{table}
\caption{Global loops for rank 4 quivers of weight 4 }
  \begin{tabular}{ l l }
    \hline
    \bf Fully cyclic classes of weight 4 & \  \  \  \bf Global mutation loops\\
    \hline
    $ \xymatrix{\cdot_{i}\ar[d]_{2}&  \ar[l]_{(2,2)} \cdot_{l} \\
\cdot_{t}\ar[ur]^{2} \ar[d]_{2} & \ar[l]_{2}\cdot_{j} \\
	\cdot_{k}  \ar[ur]_{(2. 2)}} $ & \  \  \  \ \ $\mu_{i}, \mu_{l}, \mu_{t}, \mu_{j}, \mu_{k}$\\

    $\nonumber   \xymatrix{\cdot_{k}\ar@{-}[dr]\ar[dr]^{(2, 1)}\\
 \cdot_{v}  \ar[u]^{(1, 2)}\ar[dr]_{(2, 1)}& \cdot_{j}\ar[l]_{(2,2)} \\
 & \ar[u]_{(1, 2)}\cdot_{i}} $ & \  \  \ \  \ $\mu_{ki}, \mu_{jv}$\\

 $\nonumber   \xymatrix{\cdot_{k}\ar@{-}[dr]\ar[dr]\\
 \cdot_{v}  \ar[u]\ar[dr]_{(3, 1)}& \cdot_{j}\ar[l]_{(2,2)} \\
 & \ar[u]_{(1, 3)}\cdot_{i}} $ & \  \  \ \  \ $\mu_{v}, \mu_{j}$\\
    \hline
\end{tabular}
\end{table}

``$\Rightarrow$". Suppose that $Q$ is not fully cyclic. Then $[Q]$ contains a quiver with at least one leaf vertex. Therefore, $Q$ would be one of the following cases:
\begin{equation}\label{}
\nonumber Q_{c, 1},  X_{6}, X_{7}, E^{(1, 1)}_{6}, E^{(1, 1)}_{7}, E^{(1, 1)}_{8}, \widetilde{Q}^{l}_{a, 1}, \widetilde{Q}^{l}_{a, 2},  \widetilde{Q}^{1}_{a},  \ \ \text{and} \ \ \widetilde{Q}^{l}_{a}, l\geq 2.
\end{equation}

In the following, we will cover case by case.
\begin{enumerate}
\item $Q_{c, 1}$. If $rk(Q)=4$, then $[Q]$ contains a star quiver. Then  $\mathcal{M}(Q)$ has no global loops, thanks to  Part 2 of Proposition 3.7.
  If $rk(Q)=5$ then $[Q]$ contains the following quiver

 \begin{equation}\label{}
  \nonumber Q'=\xymatrix{&\cdot_{i_{1}}&  \\
\cdot& \ar[l]\cdot_{i_{2}} \ar[u]_{2}\ar[r] & \cdot_{i_{3}}\ar@{-}[r]& \cdot}
\end{equation}
  One can see that if $\mu$ is a global loop, it must contain $\mu_{i_{2}}$. However, the vertex $i_{2}$ is a mobile vertex, i.e., one can rearrange the vertices to move $i_{2}$ to be a leaf. So, the resulting quiver would not be symmetric to $Q'$ after applying $\mu$.

  \item $E^{(1, 1)}_{6},E^{(1, 1)}_{7}, E^{(1, 1)}_{8},  X_{6}$, and $X_{7}$. One can see that in each of these cases, the mutation class contains a simply-laced quiver. Then they can not have any global loops, thanks to the proof of Part 1 of this theorem and Remark 3.11.

\item $\widetilde{Q}^{l}_{a, 1}, l\geq 1$. All of these quivers are mutationally equivalent to simply-laced quivers. Now consider the case of \begin{equation}\label{}
  \nonumber Q=\xymatrix{ \cdot_{j'}\ar[r]&   \ar[d]\cdot_{v'}   & \ar[l]\cdot_{v}\ar[r]& \cdot_{j}\ar[dl]^{4}\\
&\cdot_{k'}\ar[ul]^{4}&\cdot_{t_{1}}\ar[u]},
\end{equation}
where the weight four edges are either of valuation $(2, 2)$ or $(4, 1)$ with the apparent changes in the weights of other edges. One can see that $[Q]$ contains a simply-laced quiver. Similarly to the previous case, we can mutate $Q$ to a quiver with a simply laced 3-cycle. This process takes us to the same technique used in the previous case.
\item $\widetilde{Q}^{1}_{a, 2}, \widetilde{Q}^{1}_{a}.$ Let $Q$ be the following quiver
 \begin{equation}\label{}
  \nonumber Q'=\xymatrix{\cdot_{j'}\ar[r]^{2}& \ar[d]^{2}\cdot_{v'} \ar[r] & \cdot_{v}\ar[r]^{2}& \cdot_{j}\ar[dl]^{4}\\
&\cdot_{k'}\ar[ul]^{4}&\cdot_{k}\ar[u]^{2}},
\end{equation}
where the weight $4$ is for valuations $(2,2), (1, 4)$ or $(4,1)$ with the obvious changes. One can see that $[Q]$ contains the following quivers
\begin{equation}\label{}
  \nonumber Q':\xymatrix{\cdot_{k'}\ar[r]^{2}& \cdot_{v} \ar[r] & \cdot_{v'}\ar[r]^{2}& \cdot_{j'}\ar[dl]^{4}\\
&\cdot_{k}\ar[u]^{2}&\cdot_{j}\ar[u]^{2}}, \text{and} \end{equation}   \begin{equation}\label{}
  \nonumber Q'':\xymatrix{\cdot_{k'}\ar[r]^{2}& \cdot_{v}  & \ar[l]\cdot_{v'}\ar[d]^{2}& \ar[l]_{2}\cdot_{j'}\\
&\cdot_{k}\ar[u]^{2}&\cdot_{j}\ar[ul]_{2}}
\end{equation}
 One can see that applying any sequence of mutations that do not involve $\mu_{v}$ or $\mu_{v'}$ will not alter the edge $\xymatrix{\cdot_{v}&\ar[l] \cdot _{v'}}$. Then any sequence of mutations that do not contain $\mu_{v}$ or $\mu_{v'}$ would not add or remove edges; hence it would not cause any significant change in the underlying graph of $Q$.  Then, 
  if a sequence of mutations $\mu$ does not contain either $\mu_{v}$ or $\mu_{v'}$, it will not be a loop on $Q'$. Now assume that $\mu$ contains $\mu_{v}$ or $\mu_{v'}$, then $\mu$ must break and rebuild one or both of the 3-cycles of $Q$, which would not be a mutation loop of $Q''$, considering the positions changes of the vertices of each 3-cycle. Then $Q$ has no global loops.

  \item  $\widetilde{Q}^{l}_{a, 2}, \widetilde{Q}^{l}_{a}, l\geq 2$. One can see that the following quiver
    \begin{equation}\label{}
  \nonumber Q'=
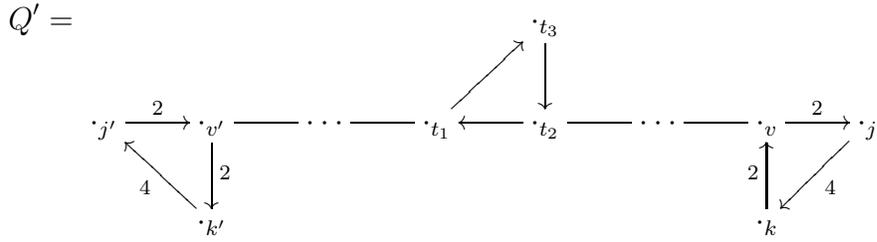
\xymatrix{&&&&\cdot_{t_{3}}\ar[d]&&&&\\
\cdot_{j'}\ar[r]^{2}&   \ar[d]^{2}\cdot_{v'} \ar@{-}[r] &\cdots\ar@{-}[r]& \cdot_{t_{1}}\ar[ur] & \ar[l]\cdot_{t_{2}} \ar@{-}[r] &\cdots &\ar@{-}[l]\cdot_{v}\ar[r]^{2} &\cdot_{j}\ar[dl]^{4}\\
&\cdot_{k'}\ar[ul]^{4}&&&&&\cdot_{k}\ar[u]^{2}}
\end{equation}
is in $[\widetilde{Q}^{l}_{a, 2}]$ and $[\widetilde{Q}^{l}_{a}]$ depending on the valuation of weight 4, where the weight 4 is of valuation either $(2, 2)$ or $(4, 1)$ with the apparent changes in the weights of other edges.  Since the subquiver $Q^*$ of $Q'$, connecting the vertices $v$ and $v'$ is simply laced, then any global loop containing one of the vertices of $Q^*$ has to contain all the vertices of $Q^*_{0}$. Where $Q^*_{0}$ is the following 3-cyclic subquiver. 
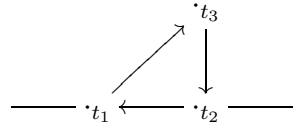
\begin{equation}\label{}
\nonumber  \xymatrix{&&\cdot_{t_{3}}\ar[d]&\\
 \ar@{-}[r]&\cdot_{t_{1}}\ar[ur] & \ar[l]\cdot_{t_{2}}\ar@{-}[r]&}
\end{equation}

 One can see that any global loop must alter the 3-cycle, $Q^*_{0}$. This argument guarantees that any global loop would contain a subsequence that breaks and rebuilds the subquiver $Q^*_{0}$. However, any simply laced quiver of three vertices has no global loop. Therefore, such a global loop $\mu$ does not exist.

\end{enumerate}
\end{enumerate}
\end{proof}
\begin{rems} \begin{enumerate}
               \item If $[Q]$ contains a simply-laced quiver then $[Q]$ does not have global loops.
               \item The case of  $rk(Q)=3$ and $w([Q])=4$ is covered in Lemma 3.3.
             \end{enumerate}

\end{rems}

 \begin{proof} The proof of both parts are  corollaries of the proof of Part 1 of Theorem 3.10
 \end{proof}

Some exceptional finite mutation types $E^{(1, 1)}_{6}, E^{(1)}_{7}, E^{(1)}_{8}$.

\begin{itemize}
  \item  [$E^{(1, 1)}_{6}$:] \begin{equation}\label{}
 \nonumber \xymatrix{ && \cdot\ar[dl]\ar[dr]\ar[drrr] &&&&\\
 \cdot\ar@{-}[r] &\cdot\ar[dr] &&\cdot\ar[dl]\ar@{-}[r]&\cdot&\cdot\ar[dlll]\ar@{-}[r]&\cdot\\
 &&\cdot\ar[uu]_{(2,2)}&&&&}
 \end{equation}

  \item  [$E^{(1, 1)}_{7}$:] \begin{equation}\label{}
 \nonumber 
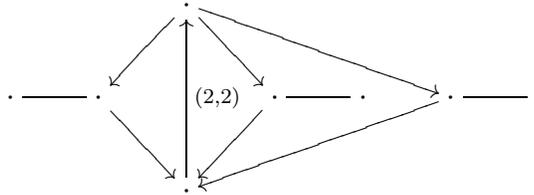
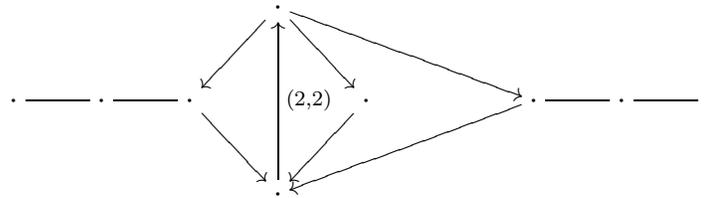
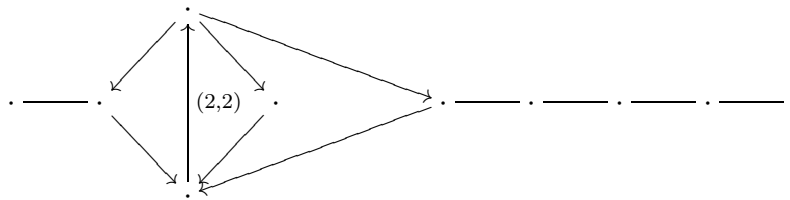
\xymatrix{ &&& \cdot\ar[dl]\ar[dr]\ar[drrr] &&&&&\\
 \cdot\ar@{-}[r]&\cdot\ar@{-}[r] &\cdot\ar[dr] &&\cdot\ar[dl]&&\cdot\ar[dlll]\ar@{-}[r]&\cdot&\ar@{-}[l]\cdot\\
 &&&\cdot\ar[uu]_{(2,2)}&&&&&}
 \end{equation}

  \item [$E^{(1, 1)}_{8}$:] \begin{equation}\label{}
 \nonumber \xymatrix{ && \cdot\ar[dl]\ar[dr]\ar[drrr] &&&&&&&\\
 \cdot\ar@{-}[r] &\cdot\ar[dr] &&\cdot\ar[dl]&&\cdot\ar[dlll]\ar@{-}[r]&\cdot&\ar@{-}[l]\cdot &\ar@{-}[l]\cdot &\ar@{-}[l]\cdot\\
 &&\cdot\ar[uu]_{(2,2)}&&&&&&&}
 \end{equation}

\end{itemize}

\subsection*{Acknowledgments}  I  would like to thank Fang Li and  Zongzhu Lin for the very valuable discussions regarding the topic of this paper.

\end{document}